\input amstex
\documentstyle{amsppt}



\define\closedconv{\operatorname{\overline{\text{conv}}}}
\define\conv{\operatorname{conv}}

\define\Un{\operatorname{Un}}
\define\ext{\operatorname{ext}}
\define\ex{\operatorname{ex}}
\define\Ch{\operatorname{Ch}}
\define\spm{\operatorname{suppmax}}

\define\Intois{\operatorname{Intois}}
\define\Is{\operatorname{Is}}

\topmatter

\title Isometries between weighted uniform spaces  \endtitle
\author Martin At. Stanev  \endauthor
\address Dept. of Mathematics, University of Chemical Technology
 and Metallurgy, 8, Kliment Ochridsky blvd., Sofia 1756, BULGARIA \endaddress
\email stanevm\@uctm.edu \endemail
\date \enddate

\abstract
The linear isometries between weighted  Banach  spaces 
of continuous functions are considered. Some of well
known theorems on isometries between spaces of continuous
functions are proved and stated, but all they are in 
an appropriate form. In this paper, we present  some
new results, too, and in  particular --- an important equivalence 
relation is investigated. 
\endabstract

\toc\widestnumber\head{A}
\head 1. Introduction\endhead
\head 2. Main Definitions and Notations\endhead
\head 3. Main Theorem\endhead
\head 4. Corollaries of the Main Theorem\endhead
\head 5. The Equivalence Relation\endhead
\head 6. Choquet Boundary\endhead
\head 7. The spaces of continuous functions $C(X)$, $C_0(Y)$\endhead
\head {}   References\endhead
\endtoc
\endtopmatter

\document

\head{1. Introduction}\endhead
We consider weighted uniform spaces and linear 
 isometries between them. A weighted uniform space 
 is a weighted Banach space of continuous functions on 
 any completely regular topological space.

In this paper a systematic studying of some of basical results on
 linear isometries between weighted uniform spaces
 is presented. Our investigation is based on theorems
 concerning the geometry of dual spaces of Banach spaces.

It is necessary to note some of previous works on linear
 isometries between spaces  classified as weghted uniform
 spaces: \cite{1}, \cite{4,~Ch.~5}  ---	onto-isometries
 of the space $C(X)$ on a Hausdorff compact; \cite{2}  ---
 onto-isometries between Banach spaces having appropriate
 function module representation; \cite{7} --- into-isometries
 of $C(X)$; \cite{9} ---  into-isometries of the space $C_0(Y)$
 of continuous functions on a locally compact Hausdorff topological
 space vanishing at infinity; \cite{14} --- linear  isometries of  the
 Bloch  spaces  of holomorphic functions in  the  unit disk;
 \cite{12} and \cite{13} --- the onto-isometries of some
 weighted Banach spaces of holomorphic functions are characterized.    
 
The main result of this paper is theorem~3.1. The basic idea
 of theorem~3.1 is demonstrated in \cite{4,~lemma~5.8.6}.
 
In section 4, corollaries of  theorem~3.1 are given. Some
 of them have analoguous assertions in the particular case
 when the corresponding Banach spaces are uniform algebras
 or are concrete Banach spaces.

In section 5, an equivalence relation is defined. Some
 restriction of this relation is used earlier in~\cite{12}.
 In our investigation, it is better to use this equivalence
 relation instead of such a relation  that is  induced 
 by a condition like that of~\cite{8,~p.195,~(2)}.

In section 6, an extension of the classical notion Choquet
 boundary is given. The extension of this notion is based
 on the key property that such a boundary characterizes the
 extreme functionals.

In section 7, we apply the above developped theory to prove
 the well known  classical results on isometries and on
 extreme functionals of spaces of continuous functions.
 We obtain the result \cite{9,~theorem~1}, too. In 
 this section we state some results that are important
 for the general theory of isometries  between weighted 
 uniform spaces, but not  only for the special case under 
 consideration in this section. 

Moreover, it is possible to obtain theorems
 on isometries such as in~\cite{1}, \cite{7}, \cite{9} and
 theorems on extreme functionals as in~\cite{4,~lemma~5.8.6,
 lemma~5.8.7}, by using the presented here results. 

We consider in this paper spaces of complex-valued functions.
 The proofs are given with the set $S$ of all complex numbers
 of modulus one. Exactly the same results are valide for the
 spaces of real-valued functions, but in this case the set
 $S=\{-1,1\}$ is used instead of the set of all complex numbers
 of modulus one. The proofs are the same in both cases.

\remark{Remark}Some notations and some terminology explanations
 are placed just before their use. In each section it is used
 the terminology, notations and results of the previous sections.
\endremark

\head{2. Main Definitions and Notations} \endhead
We use the standart terminology and notations of the theory
 of Banach spaces but it is necessary to give some explanations.
 
According to~\cite{11,~\S~ 3.20}, the definitions of extreme
 subsets and extreme points are next. Let~$L$ be any complex
 linear space. Suppose~$E$ is a non empty subset of~$L$. We recall
 that a subset~$D$, $\emptyset\ne D\subset E$, is said to be
 an extreme subset of~$E$ iff it is fulfilled the following condition:
\block
 for every $x,y\in E$ such that there is a number~$t$, $0<t<1$,
 for which $tx+(1-t)y\in D$, in fact it is true that $x,y\in D$.
\endblock
 A point $x_0$ of~$E$ is an extreme point of~$E$ iff the
 single point set $\{x_0\}$ is an extreme subset of~$E$.
 
We make use of notations which are not in common use.

The set of all extreme points of~$E$ is denoted by~$\ex E$.
 The convex hull of a set~$E\subset L$ is denoted by~$\conv E$.
 If, in addition,~$L$ is a Hausdorff topological linear space,
 then  the closed convex hull of~$E$ is $\closedconv E$.

We recall that if~$E$ is a compact subset of a Hausdorff
 topological linear space then:
\roster
\item"{$\bullet$}" $\ex E\ne \emptyset$;
\item"{$\bullet$}" $E\subset \closedconv (\ex E)$ and,
 in particular $E=\closedconv (\ex E)$ when $E$ is convex;
\item"{$\bullet$}" $\ex (\closedconv E)\subset E$, when
 $\closedconv E$ is compact.
\endroster
 Note that the first two bullets represent the Krein-Milman
 theorem and the third bullet is a well-known assertion
 (see for instance one of the following two monographs   
 \cite{5,~proposition~10.1.3} or \cite{4,~lemma~5.8.5}).

Let $(A,\|\cdot\|_A)$ be an abstract complex Banach space with
 the norm~$\|\cdot\|_A$. Then $(A^*,\|\cdot\|_{A}^*)$ is
 the corresponding dual Banach space. In this paper the
 following specifical notations are used:
$$
\align
 (A)_1&=\bigl\{\,f:f\in A,\;\|f\|\le1\bigr\},\\
\partial (A)_1&=\bigl\{\,f:f\in A,\;\|f\|=1\bigr\},
\endalign
$$
 and for brevity it is setted $\ext A^*=\ex (A^*)_1$.

We use some notations on isometries which are not in
 common use. Let for $i=1,2$, the abstract Banach
 space $(A_i,\|\cdot\|_{A_i})$ be with the norm
 $\|\cdot\|_{A_i}$. Then we write
$$
 T\in\Intois(A_1,A_2)
$$
 iff the linear operator $T:A_1\rightarrow A_2$ is such
 that $\|Tf\|_{A_2}=\|f\|_{A_1}$, $\forall f\in A_1$.
 We set $\Is(A_1,A_2)$ to be the family of all
 $T\in\Intois(A_1,A_2)$ which maps $A_1$ onto $A_2$.
 It follows by the Krein-Milman theorem that for the
 dual map $T^*$ of $T$, $T\in\Intois(A_1,A_2)$, we have
$$
 T^*(\ext A_2^*)\supset \ext A_1^*,
$$
 i.e., $\forall\ell\in\ext A_1^*$ there exists $m\in\ext A_2^*$
 such that $T^*(m)=\ell$.

Note that if $F:\Omega_1\rightarrow\Omega_2$ is an
 abstract map from $\Omega_1$ to $\Omega_2$ then as usual one set
$$
 F(Y)=\{F(y)|\; y\in Y\},
$$
 where $\emptyset\ne Y\subset\Omega_1$, and
 $F(\emptyset)=\emptyset$.

The set of all scalars of modulus one is denoted by $S$,
 thus in the complex case
$$
S=\{\,\lambda:\lambda\in\Bbb C,\;|\lambda|=1\}
$$
(and in the real case $S=\{1,-1\}$).

\definition{Definition 2.1}We define that the triple
 $(Z,\Delta_A,A)$ has the property $(*)$ if and only
 if $Z$ is a set, $A$ is a non zero Banach space,
 $\Delta_A:S\times Z\rightarrow (A^*)_1$ is a single-valued
 map such that:
\roster
\item $\displaystyle\|f\|_A=\sup_{z\in Z}|\Delta_A(1,z)(f)|,
 \quad\forall f\in A$;
\item $\Delta_A(\lambda,z)=\lambda\Delta_A(1,z),
 \quad\forall\lambda\in S,\forall z\in Z$.
\endroster
\enddefinition

Let the triple $(Z,\Delta_A,A)$ have the property $(*)$.
 We shall say that $\Delta_A$ is continuous when $Z$ is
 a topological space, $(A^*)_1$ has the induced
 weak${}^*$-topology of $A^*$ and the map
 $\Delta_A:S\times Z\rightarrow (A^*)_1$ is continuous.

Let the triple $(Z,\Delta_A,A)$ have the property $(*)$.
 We shell use the following version of the notion
 ``boundary of the space A'': a set~$Y$,~$Y\subset Z$,
 is said to be a boundary of~$Z$ with respect to the
 space~$A$ (or briefly, $Y$ is a boundary of~$A$) iff
$$
\|f\|_A=\sup_{y\in Y}|\Delta_A(1,y)(f)|,\qquad \forall f\in A.
$$

An example of such a triple $(Z,\Delta_A,A)$ that has
 the property $(*)$ is constructed as follows. Let~$Z$
 be any completely regular topological space and let
 $p:Z\rightarrow\Bbb C\setminus\{0\}$ be any continuous
 function. The Banach space $(C(Z;p),\|\cdot\|_p)$
 is defined by
$$
 C(Z;p)=\{f|\;f:Z\rightarrow\Bbb C\text{ is continuous, }
 \|f\|_p:=\sup_{z\in Z}|p(z)f(z)|<\infty\}.
$$
 The notation $\Un(Z;p)$ stands for the family of all
 non zero Banach subspaces of~$C(Z;p)$. We say that the
 Banach space~$A$ is weighted uniform space
 iff~$A\in\Un(Z;p)$, where $Z$ is any completely regular
 topological space and $p:Z\rightarrow\Bbb C\setminus\{0\}$
 is any continuous function.
 If~$A\in\Un(Z;p)$ then the map
$$
 \Delta_A:S\times Z\rightarrow (A^*)_1
$$
 is defined by $\Delta_A(\lambda,z)(f)=\lambda p(z)f(z)$,
 where $\lambda\in S$, $z\in Z$, $f\in A$. Then the triple
$(Z,\Delta_A,A)$ has the property $(*)$ and moreover
 $\Delta_A$ is continuous.

\head{3. Main Theorem} \endhead
In this section we state and prove our main theorem---
 theorem 3.1. In its proof it is used one of the ideas
 which are demonstrated in \cite{4,~lemma~5.8.6}.

\proclaim{Theorem 3.1}Let $A_1$ be a non zero Banach
 space and let the triple $(Z_2,\Delta_{A_2},A_2)$ have
 the property $(*)$. Suppose $\Intois(A_1,A_2)\ne\emptyset$.
 If~$T\in\Intois(A_1,A_2)$ then
$$
 \ext A_1^*\subset T^*(\overline{\;\Delta_{A_2}
 (S\times Y)\;}),
$$
 where $Y$ is any boundary of $A_2$ and
 $\overline{\Delta_{A_2}(S\times Y)}$ is the closure
 of~$\Delta_{A_2}(S\times Y)$ in~$A_2^*$ with respect
 to the weak${}^*$-topology.
\endproclaim
\demo{Proof of theorem 3.1}Note, $\ext A_1^*\ne\emptyset$
 by the Krein-Milman
 theorem.

Let $T\in\Intois(A_1,A_2)$ be any and let $Y$ be any
 boundary of the space~$A_2$. Then
 $\forall\ell\notin\closedconv T^*(\Delta_{A_2}
 (S\times Y))$ there is an element $f\in A_1$ such that
$$
 Re\;\ell(f)>\sup_{m\in\closedconv T^*(\Delta_{A_2}
 (S\times Y))}Re\;m(f).
$$
 Hence
$$
\align
 Re\;\ell(f)>\sup_{m\in T^*(\Delta_{A_2}(S\times Y))}
 Re\;m(f)&=\sup_{(\lambda,z)\in S\times Y} Re\;
 \Delta_{A_2}(\lambda,z)(Tf)\\
 {}&=\|Tf\|_{A_2}=\|f\|_{A_1}.
\endalign
$$
 Therefore $\ell\notin (A_1^*)_1$ and we obtain
$$
 (A_1^*)_1\subset\closedconv T^*(\Delta_{A_2}(S\times Y)).
$$
 The inclusion
$$
 (A_1^*)_1\supset\closedconv T^*(\overline{\;
 \Delta_{A_2}(S\times Y)\;}).
$$
 follows by an obvious way from $T\in\Intois(A_1,A_2)$
 and we omit the details. Thus
$$
 (A_1^*)_1=\closedconv T^*(\overline{\;\Delta_{A_2}
 (S\times Y)\;}).
$$
 Then by the weak${}^*$-compactness of~$T^*(\overline{\;
 \Delta_{A_2}(S\times Y)\;})$ it follows
$$
 \ext A_1^*\subset T^*(\overline{\;\Delta_{A_2}
 (S\times Y)\;}).
$$
 So, theorem 3.1 is proved.
\enddemo

\remark{Remark}Suppose the triple $(Z,\Delta_A,A)$ has
 the property $(*)$ and in addition, let the map $\Delta_A$
 be contiuous. If $K$ is a compact
 subset of $Z$ then, of course, the image
 $\Delta_A(S\times K)$ is weak${}^*$-compact. So that
$$
 \overline{\Delta_{A}(S\times K)}=\Delta_{A}(S\times K),
$$
 where $\overline{\Delta_{A}(S\times K)}$ is the
 weak${}^*$-closure of $\Delta_{A}(S\times K)$ in the
 dual space $A^*$.
\endremark

\head{4. Corollaries of the Main Theorem} \endhead
In this section are placed corollaries of theorem 3.1
 which are important for further results.

It is necessary to explique the following notations.
 Let $A$ be any non zero Banach space. For every one non
 empty family of its elements $G$, $\emptyset\ne G\subset
 A\setminus\{0\}$, we set
$$
 \Sigma^A(G)=\{\ell|\;\ell\in\partial(A^*)_1,
 |\ell(f)|=\|f\|_A\text{ for every }f\in G\/\}.
$$
 We define that the family $G$ is centered iff
 $\Sigma^A(G)\ne\emptyset$.

Let the triple $(Z,\Delta_A,A)$ have the property $(*)$.
 We introduce the notation $\spm(f)$ as follows:
$$
 \spm(f)=\{z|\;z\in Z,|\Delta_A(1,z)(f)|=\|f\|_A\},
 \qquad f\in A.
$$
 We define that the family $G$, $\emptyset\ne G\subset
 A\setminus\{0\}$, is placed over $V$, where $\emptyset\ne
 V\subset Z$, iff there is an element $f\in G$ such that
$$
 \sup_{z\in Z\setminus V}|\Delta_A(1,z)(f)|<\|f\|_A.
$$

\proclaim{Lemma 4.1}Let $E$ be any non empty subset
 of a complex linear space. If $D$ is an extreme subset
 of $E$, then
$$
 \ex D=D\cap\ex E.
$$
\endproclaim

Instead of a detailed proof, we have to remark that
 lemma~4.1 is an immediate consequence of the definition
 of extreme subsets.

The inclusion $\ex D\subset D\cap\ex E$ follows from
 the fact that $D$ satisfies the condition characterizing
 extreme subsets (indeed, for such a $D$ every one of its
 extreme points is an extreme point of $E$, too);
 the case $\ex D=\emptyset$ is a trivial one. The opposite
 inclusion $\ex D\supset D\cap\ex E$ is an obvious
 consequence of $D\subset E$.

\proclaim{Lemma 4.2}Let the triple $(Z,\Delta_A,A)$ have
 the property $(*)$. If the family $G$, $\emptyset\ne
 G\subset A\setminus\{0\}$, is centered, then
\roster
\item $\ex\Sigma^A(G)\ne\emptyset$;
\item $\Sigma^A(G)\subset\closedconv(\ex\Sigma^A(G))$;
\item $\ex\Sigma^A(G)=\Sigma^A(G)\cap\ext A^*$.
\endroster
\endproclaim
\demo{Proof of lemma 4.2}Let $G$, $\emptyset\ne
 G\subset A\setminus\{0\}$, be centered. Then,
 in particular,
$$
 \emptyset\ne\Sigma^A(G)\subset (A^*)_1
$$
 and $\Sigma^A(G)$ is weak${}^*$-closed. Therefore
 $\Sigma^A(G)$ is weak${}^*$-compact and by the
 Krein-Milman theorem we obtain assertions
 \therosteritem1 and \therosteritem2 of this lemma. In addition,
 by a direct computation, we verify that $\Sigma^A(G)$
 is an extreme subset of $(A^*)_1$. Then the assertion
 of lemma~4.2.\therosteritem3 is the particular case
 of lemma~4.1 when $E=(A^*)_1$, $D=\Sigma^A(G)$.

Thus lemma~4.2 is proved.
\enddemo

\proclaim{Theorem 4.1}Let $A_1$ be a non zero Banach
 space, the triple $(Z_2,\Delta_{A_2},A_2)$
 have the property $(*)$ and let $\Delta_{A_2}$ be a
 continuous map. Suppose, in addition,
 $\Intois(A_1,A_2) \ne\emptyset$ and let 
$T\in\Intois(A_1,A_2)$ be any.
 If the centered family $G$, $\emptyset\ne G\subset A_1
 \setminus\{0\}$, is such that $TG$ is placed over
 some compact subset of $Z_2$ then
$$
\ex\Sigma^{A_1}(G)\subset T^*\Delta_{A_2}(S\times
 (\bigcap_{f\in G}\spm(Tf))) \bigcap\ext A_1^*,
$$
and in particular, there exists $\displaystyle z_0\in
 \bigcap_{f\in G}\spm(Tf)$ for which we have
$$
T^*\Delta_{A_2}(1,z_0)\in\ext A_1^*.
$$
\endproclaim

\demo{Proof of theorem 4.1}Suppose  the family $G$,
 $\emptyset\ne G\subset A_1\setminus\{0\}$, is 
 centered and  in addition,  it is such that $TG$ is placed over
 some compact subset of $Z_2$.  Let $K$ stand for such a compact set. Then
 there is $f_0\in G$ such that
$$
 \sup_{z\in Z_2\setminus K}|\Delta_{A_2}(1,z)
 (Tf_0)|<\|Tf_0\|_{p_2}.
$$
 Thus $f_0$ represents a weak${}^*$-continuous
 linear functional defined on $A_1{}^*$ such that
 distinguishes $\Sigma^{A_1}(G)$ and
 $T^*(\overline{\;\Delta_{A_2}(S\times
 (Z_2\setminus K))\;})$, (here, the overlined set stands 
 for the weak${}^*$-\-closure of $\Delta_{A_2}(S\times
 (Z_2\setminus K))$ in the space $A_2{}^*$). So,
$$
 \Sigma^{A_1}(G)\cap T^*(\overline{\;\Delta_{A_2}
 (S\times (Z_2\setminus K))\;})=\emptyset.\tag{4.1}
$$
Moreover, by lemma 4.2 when it is setted $Z=Z_1$, $p=p_1$,
 $A=A_1$, we obtain for $G$
$$
 \emptyset\ne\ex\Sigma^{A_1}(G)=
 \Sigma^{A_1}(G)\cap\ext A_1{}^*.
$$
Further, it follows by theorem 3.1
$$
 \ex\Sigma^{A_1}(G)\subset\Sigma^{A_1}(G)
 \cap T^*(\overline{\;\Delta_{A_2}(S\times Z_2)\;})\cap\ext A_1{}^*,
$$
 and hence
 $$
 \ex\Sigma^{A_1}(G)\subset\Sigma^{A_1}
 (G)\cap (T^*\Delta_{A_2}(S\times K)\cup
 T^*(\overline{\,\Delta_{A_2}(S\times
 (Z_2\setminus K))\,}))\cap\ext A_1{}^*.
$$
Then, by the equation \thetag{4.1}
$$
 \ex\Sigma^{A_1}(G)\subset\Sigma^{A_1}(G)\cap
 T^*(\Delta_{A_2}(S\times K))\cap\ext A_1{}^*.
$$
 and hence
$$
  \ex\Sigma^{A_1}(G)\subset T^*(\Delta_{A_2}(S\times
 (K\cap (\bigcap_{f\in G}\spm(Tf)))))\bigcap\ext A_1{}^*.
$$
In particular, the asserted inclusion in theorem~4.1 is  proved. 
\enddemo

\proclaim{Corollary 4.1}Let $A_1$ be a non zero Banach
 space and let the triple $(Z_2,\Delta_{A_2},A_2)$
 have the property $(*)$. Let $Z_2$ be a Hausdorff
 topological space and the map $\Delta_{A_2}$ be
 continuous with respect to the weak${}^*$-topology
 of $A_2^*$. Suppose, in addition, $\Intois(A_1,A_2)
 \ne\emptyset$ and let $T\in\Intois(A_1,A_2)$ be any.
 If the family $G$, $\emptyset\ne G\subset A_1
\setminus\{0\}$, is such that
\roster
\item"{(i)}" $\displaystyle\bigcap_{f\in G}
 \spm(Tf)=\{z_0\}$, where $z_0$ is a point of $Z_2$,
\item"{(ii)}" the image $TG=\{Tf|\;f\in G\/\}$ of
 $G$, is placed over a compact subset of $Z_2$,
\endroster
 then
$$
 \Sigma^A(G)=T^*\Delta_{A_2}(S\times\{z_0\})
 \subset\ext A_1{}^*.
$$
\endproclaim
\demo{Proof of corollary 4.1}Let for the family $G$,
 $\emptyset\ne G\subset A_1\setminus\{0\}$, be
 fulfilled the conditions~\therosteritem{i}
 and~\therosteritem{ii}. In particular, from the
 condition~\therosteritem{i} it follows
 $T^*\Delta_{A_2}(1,z_0)\in\Sigma^{A_1}(G)$,
 so that $\Sigma^{A_1}(G)\ne\emptyset$. Thus the
 family $G$ is centered. 
Then by theorem~4.1 it follows 
 $\ex\Sigma^{A_1}(G)\subset
 T^*(\Delta_{A_2}(S\times\{z_0\}))\cap\ext  A_1^*$, and hence
$$
 \ex\Sigma^{A_1}(G)=T^*(\Delta_{A_2}(S\times\{z_0\}))\subset \ext A_1^*.
$$
  Further, from the Krein-Milman theorem
$$
\align
\Sigma^{A_1}(G)&\subset\closedconv
 T^*(\Delta_{A_2}(S\times\{z_0\}))=\\
	     &=\{rT^*(\Delta_{A_2}(\lambda,z_0))|\;
	      0\le r\le 1,\lambda\in S\/\}.\tag{4.2}
\endalign
$$
Moreover, accordingly to the definition of
 $\Sigma^{A_1}(G)$, the equation
 $|\ell(f)|=\|f\|_{A_1}$ holds for every
 $\ell\in\Sigma^{A_1}(G)$ and $\forall f\in G$.
 Then by the inclusion \thetag{4.2} we obtain
$$
 \Sigma^{A_1}(G)=T^*(\Delta_{A_2}(S\times\{z_0\})).
$$
 Thus corollary 4.1 is proved.
\enddemo

\proclaim{Proposition 4.1}Let $A_1$ be a non zero Banach
 space and let the triple $(Z_2,\Delta_{A_2},A_2)$
 have the property $(*)$. Suppose, in addition,
 $\Intois(A_1,A_2)\ne\emptyset$ and let
 $T\in\Intois(A_1,A_2)$ be any. If the set of
 functionals $\Cal L$, $\emptyset\ne\Cal L\subset
 (A_1^*)_1$, is such that there are $f,Q$,  $f\in A_1$
 and $Q$ is a subset of $Z_2$, such that at least one of
 the following two inequalities holds:
\roster
\item"{$\bullet$}" $\displaystyle
 \sup_{z\in Z_2\setminus Q}|T^*
 (\Delta_{A_2}(1,z))(f)|<
 \inf_{\ell\in\Cal L}|\ell(f)|$,
\item"{$\bullet$}" $\displaystyle
 \inf_{z\in Z_2\setminus Q}|T^*
 (\Delta_{A_2}(1,z))(f)|>
 \sup_{\ell\in\Cal L}|\ell(f)|$,
\endroster
 then $\overline{\Cal L}\cap T^*(\overline{\;
 (\Delta_{A_2}(S\times (Z_2\setminus Q)))\;})=\emptyset$.
\endproclaim
 The assertion is an immediate consequence of the
 assumptions and we omit its proof.

\proclaim{Corollary 4.2}Let $A_1$ be a non zero Banach
 space and let the triple $(Z_2,\Delta_{A_2},A_2)$ have
 the property $(*)$. Let $Z_2$ be a Hausdorff topological
 space and the map $\Delta_{A_2}$ be continuous with
 respect to the weak${}^*$-topology of $A_2^*$. Suppose,
 in addition, $\Intois(A_1,A_2)\ne\emptyset$ and let
 $T\in\Intois(A_1,A_2)$ be any. If the point
 $z_0\in Z_1$ is such that $\Delta_{A_1}(1,z_0)\in
 \ext A_1{}^*$ and if  there are $f,K$, $f\in
 A_1\setminus\{0\}$ and $K$ is a compact subset of
 $Z_2$, for which at least one of the following two
 inequalities holds:
\roster
\item"{$\bullet$}" $\dsize\sup_{z\in Z_2\setminus
 K}|T^*(\Delta_{A_2}(1,z)(f))|<|\Delta_{A_1}(1,z_0)(f)|$,
\item"{$\bullet$}" $\dsize\inf_{z\in Z_2\setminus
 K}|T^*(\Delta_{A_2}(1,z)(f))|>|\Delta_{A_1}(1,z_0)(f)|$,
\endroster
then
\roster
\item $\Delta_{A_1}(S\times\{z_0\})\cap T^*(\overline{\;
 (\Delta_{A_2}(S\times (Z_2\setminus K)))\;})=\emptyset$,
\item $\Delta_{A_1}(S\times\{z_0\})\subset
 T^*(\Delta_{A_2}(S\times K))$.
\endroster
\endproclaim
 Instead of a proof, note that corollary~4.2 is an
 immediate consequence of both, theorem~3.1 and
 proposition~4.1.

\head{5. The Equivalence Relation}\endhead
 In this section, it is defined the equivalence
 relation $\thicksim_A$ between the points of the
 set $Z$, where $(Z,\Delta_A,A)$ is any triple having
 the property $(*)$. This relation is induced by the
 map $\Delta_A:S\times Z\rightarrow A^*$. It is
 important to note that there is triple $(Z,\Delta_A,A)$
 such that between its   points $x,y\in Z$ we have
 $x\thicksim_A y$ iff the functions of $A$ do not
 distinguish points $x$ and $y$. For instance, we
 have such an equivalence when $Z$ is any completely
 regular topological space and $A\in \Un(Z;1)$ is
 such that $A\ni 1$.

\definition{Definition 5.1}Let the triple
 $(Z,\Delta_A,A)$ has the property $(*)$. The relation
 $\thicksim_A$ between points of $Z$ is defined by
$$
 x\thicksim_A y\Longleftrightarrow|
 \Delta_A(1,x)(f)|=|\Delta_A(1,y)(f)|,\quad\forall
 f\in A.
$$
\enddefinition

\remark{Remark}A direct computation shows that the
 above defined relation is an equivalence relation.
\endremark

\definition{Definition 5.2}Let the triple
 $(Z,\Delta_A,A)$ has the property $(*)$.
 Suppose $V$ is any non empty subset of $Z$.
 We define that $|pA|$ distinguishes  the points
 of $V$ iff for every couple $x,y\in V$, $x\ne y$,
 there is an element $f\in A$ such that
$$
 |\Delta_A(1,x)(f)|\ne |\Delta_A(1,y)(f)|.
$$
\enddefinition

\proclaim{Proposition 5.1}Let the triple
 $(Z,\Delta_A,A)$ has the property $(*)$. Suppose
 $x,y\in Z$. Then
$$
x\thicksim_A y \Longleftrightarrow\text{ there are
 }\lambda,\mu\in S\text{ such 
that }\Delta_A(\lambda,x)=\Delta_A(\mu,y).
$$
\endproclaim
\demo{Proof of proposition 5.1}
 Note, the assertion ``$\Longleftarrow$''
 is trivial and we omit its proof.

Proof of the part ``$\Longrightarrow$''.
 Thus we assume $x\thicksim_A y$. It is necessary
 to establish that there are scalars
 $\lambda,\mu\in S$ such that
$$
 \Delta_A(\lambda,x)=\Delta_A(\mu,y).\tag{5.1}
$$
 According to the definition 5.1, it follows
 from $x\thicksim_A y$ that
$$
 |\Delta_A(1,x)(f)|=|\Delta_A(1,y)(f)|,
 \quad\forall f\in A.\tag{5.2}
$$
 Hence for the kernels $\ker\Delta_A(1,x)$ and
 $\ker\Delta_A(1,y)$ of functionals $\Delta_A(1,x)$
 and $\Delta_A(1,y)$ we have
$$
 \ker\Delta_A(1,x)=\ker\Delta_A(1,y).
$$
 Therefore, there are two scalars $z,w\in\Bbb C$,
 at least one of them is not zero, and such that
$$
 z\Delta_A(1,x)+w\Delta_A(1,y)=0.\tag{5.3}
$$

Further, we consider the following two hypothesises
 about the functional $\Delta_A(1,y)$:
 $\Delta_A(1,y)=0$, and $\Delta_A(1,y)\ne 0$.

First, suppose $\Delta_A(1,y)=0$. Then
 $\Delta_A(1,y)(f)=0$, $\forall f\in A$. By the
 equation \thetag{5.2} it follows $\Delta_A(1,x)(f)=0$,
 $\forall f\in A$, i.e., $\Delta_A(1,x)=0$. Hence
 $\Delta_A(1,x)=\Delta_A(1,y)$ and the equation
 \thetag{5.1} is established with $\lambda=\mu=1$.

Second, suppose $\Delta_A(1,y)\ne 0$. Then,
 in particular, there is an element $f_0\in A$
 such that $\Delta_A(1,y)(f_0)\ne 0$. By the
 equation \thetag{5.2} it follows
 $\Delta_A(1,x)(f_0)\ne 0$. In addition, according
 to the equation \thetag{5.3} we obtain
$$
 z\Delta_A(1,x)(f_0)+w\Delta_A(1,y)(f_0)=0.
$$
Then from the equation \thetag{5.2} it follows
 $|z|=|w|$ and in particular $z\ne 0$. Thus
 $wz^{-1}\in S$. Then by \thetag{5.3}
$$
 \Delta_A(1,x)=\Delta_A(-wz^{-1},y),
$$
 i.e., the equation \thetag{5.1} is established
 with $\lambda=1$, $\mu=-wz^{-1}$. So, the
 assertion ``$\Longrightarrow$'' is proved and
 proposition 5.1 is proved, too.
\enddemo

\proclaim{Proposition 5.2}Let the triple
 $(Z,\Delta_A,A)$ has the property $(*)$.
 Suppose $L$ and $M$ are subsets of $Z$ such
 that $\emptyset\ne L\subset M\subset Z$. Then
 $\forall x\in L$, $\forall y\in M$, $x\ne y$,
 we have $x\not\thicksim_A y$ if and only if the
 following two conditions are fulfilled:
\roster
\item the restriction $\Delta_A:S\times
 L\rightarrow (A^*)_1$ is an injective map;
\item $\Delta_A(S\times(M\setminus L))\cap
 \Delta_A(S\times L)=\emptyset$.
\endroster
\endproclaim
 Instead of a detailed proof we have to
 remark that proposition~5.2 is a direct
 consequence of proposition~5.1.

\definition{Definition 5.3}Let the triple
 $(Z,\Delta_A,A)$ has the property $(*)$. Then
 we define that the couple of sets $(M,L)$,
 $\emptyset\ne L\subset M\subset Z$, is proper
 with respect to $(Z,\Delta_A,A)$ iff the following
 three conditions are fulfilled:
\roster
\item the restriction $\Delta_A:S\times
 L\rightarrow (A^*)_1$ is an injective map;
\item the restriction $\Delta_A:S\times
 M\rightarrow \Delta_A(S\times M)$ is closed,
 where $\Delta_A(S\times M)$ has the induced
 weak${}^*$-topology of $A^*$;
\item $\Delta_A(S\times(M\setminus L))\cap
 \Delta_A(S\times L)=\emptyset$.
\endroster
\enddefinition

\remark{Remark 5.1}Let the triple $(Z,\Delta_A,A)$
 has the property $(*)$. Suppose the couple of sets
 $(M,L)$, $\emptyset\ne L\subset M\subset Z$, is
 such that $\forall x\in L$, $\forall y\in M$,
 $x\ne y$, we have $x\not\thicksim_A y$ and in
 addition, let the restriction
$$
 \Delta_A:S\times M\rightarrow
 \Delta_A(S\times M)\tag{5.4}
$$
 be closed with respect to the induced
 weak${}^*$-topology of $A^*$ . Then the restriction
$$
 \Delta_A:S\times L\rightarrow
 \Delta_A(S\times L)\tag{5.5}
$$
 is a homeomorphic map. The following arguments
 prove this assertion. By proposition~5.2, the
 restriction \thetag{5.5} is bijective. Moreover,
 by proposition~5.2, it is the restriction on the
 whole preimage of $\Delta_A(S\times L)$ under the
 map \thetag{5.4} and hence it is closed, too.
 Then it is a homeomorphic map.
\endremark

\head{6. Choquet Boundary}\endhead
In this section, an extension of the classical
 notion Choquet boundary is given.

\definition{Definition 6.1}Let $A_1$ be a non
 zero Banach space and let the triple $(Z_2,
 \Delta_{A_2},A_2)$ have the property $(*)$.
 Suppose $\Intois(A_1,A_2)\ne\emptyset$ and
 let $T\in\Intois(A_1,A_2)$ be any. We define
 the following two notations:
\roster
\item"{$\bullet$}" $M_T(A_1)=\{z|\;z\in Z_2,
 T^*\Delta_{A_2}(1,z)\in\ext A_1{}^*\}$;
\item"{$\bullet$}" $\Ch_T(A_1)$ is the family
 of all non empty subsets $Y$ of $Z_2$ such
 that the image $T^*\Delta_{A_2}
 (S\times Y)=\ext A_1{}^*$.
\endroster
 Moreover, instead of $M_T(A_1)$ and of
 $\Ch_T(A_1)$ it is written $M(A)$ and $\Ch(A)$
 resp., when $A_1=A_2=A$, and $T:A\rightarrow A$
 is the identity. In the more special case when
 $Z$ is a Hausdorff compact, $A\in\Un(Z;1)$ and
 $A\ni 1$, one define that a subset $Y$ of $Z$
 is said to be Choquet boundary iff $Y\in\Ch(A)$.
\enddefinition

\remark{Remark}If $\Ch_T(A_1)\ne\emptyset$ then
 it is obvious that $M_T(A_1)\in\Ch_T(A_1)$.
\endremark

\remark{Remark 6.1}Suppose $A_1$ is a non
 zero Banach space and let the triple $(Z_2,
 \Delta_{A_2},A_2)$ have the property $(*)$.
 If  $T\in\Intois(A_1,A_2)$ is such that $\Ch_T(A_1)\ne\emptyset$.  Then:
\roster
\item $\Ch_T(A_1)=\Ch(TA_1)$;
\item $M_T(A_1)=M(TA_1)$. 
\endroster
We sketch the proof by the following helping arguments. 
For brevity let's put $I:TA_1\to A_2$ to be the identity embedding operator, let's define $\Delta_{TA_1}=I^*\Delta_{A_2}$,  and let $\tilde T:A_1\to TA_1$ be the restriction of the isometry $T\in\Intois(A_1,A_2)$. In particular, $\tilde T\in\Is(A_1,TA_1)$ and by the Krein--Milman theorem ${\tilde T}^*:\ext (TA_1)^*\to \ext A_1^*$ is bijective. Further by the factorization $T^*={\tilde T}^*I^*$ it follows 
$$
T^*\Delta_{A_2}={\tilde T}^*I^*\Delta_{A_2}={\tilde T}^*\Delta_{TA_1}. 
$$
A direct consequence of this equation are both \therosteritem1 and \therosteritem2, and we omit the details.
\endremark

\proclaim{Proposition 6.1}Let $A_1$ be a non
 zero Banach space and let the triple $(Z_2,
 \Delta_{A_2},A_2)$ have the property $(*)$.
 Then the following two assertions hold:
\roster
\item in the case when $\Delta_{A_2}$ is continuous and $A_2$ 
has a compact boundary $K\subset Z_2$, if  $T\in\Intois(A_1,A_2)$ ,then
$$
\ext A_1{}^*\subset T^*\Delta_{A_2}(S\times K),
$$
 and hence, in particular, $M_T(A_1)\cap K\ne
 \emptyset$ and
$$
M_T(A_1)\cap K\in\Ch_T(A_1);
$$
\item in the case when  $\Ch(A_2)\ne\emptyset$ and $B$ 
is a non zero Banach subspace of $A_2$, if $T\in\Intois(A_1,B)$, then
for every $Y\in\Ch(A_2)$  (below, $\Delta_B(S\times Y)$ stands for the set of restrictions of functionals of $\Delta_{A_2}(S\times Y)$ on the subspace $B$ of $A_2$)
$$
 \ext A_1{}^*\subset T^*\Delta_B(S\times Y),
$$
 and hence, in particular, $M_T(A_1)\cap
 Y\ne\emptyset$ and
$$
M_T(A_1)\cap Y\in\Ch_T(A_1).
$$
\endroster
\endproclaim

\demo{Proof of proposition 6.1}Instead of an
 explicite proof of assertion \therosteritem1,
 we have to remark that it is an immediate
 consequence of theorem~3.1 and we omit the
 details.

Proof of assertion \therosteritem2. Let $B$ 
be a non zero Banach subspace of $A_2$. Let the
 operator $I:B\rightarrow A_2$ stand
for the identity embedding of $B$ into  $A_2$. 
Then in particular, the triple $(Z_2,\Delta_B,B)$ has 
the property $(*)$, where $\Delta_B$ is defined by 
$\Delta_B=I^*\Delta_{A_2}$.  Note $I\in\Intois(B,A_2)$ 
 so, that  $\ext B^*\subset I^*(\ext  A_2{}^*)$.
Then
$$
\align
\ext A_1{}^*&\subset T^*(\ext B^*)\subset
 T^*(I^*(\ext A_2{}^*))\\
		    &\subset T^*(I^*(\Delta_{A_2{}^*}
 (S\times Y)))=T^*(\Delta_B(S\times Y)),
\endalign
$$
and so, assertion \therosteritem2 is proved.

Thus proposition 6.1 is proved, too.
\enddemo

\proclaim{Proposition 6.2}Let for $i=1,2$ the
 triples $(Z_i,\Delta_{A_i},A_i)$ have the
 property $(*)$. Suppose that $\Ch(A_1)\ne
\emptyset$. Suppose, in addition,
$\Intois(A_1,A_2)\ne\emptyset$ and let
 $T\in\Intois(A_1,A_2)$ be any.
If the set $V$, $\emptyset\ne V\subset
 Z_1$, is such that there are $f,W$,
  $f\in A_1$ and $W$ is a subset of $Z_2$,
 such that at least one of the following two
 inequalities holds:
\roster
\item"{$\bullet$}" $\displaystyle
\sup_{z\in W} |T^*(\Delta_{A_2}(1,z))(f)|<
\inf_{x\in Z_1\setminus V}|\Delta_{A_1}(1,x)(f)|$,
\item"{$\bullet$}" $\displaystyle
\inf_{z\in W }|T^*(\Delta_{A_2}(1,z))(f)|>
\sup_{x\in Z_1\setminus V}|\Delta_{A_1}(1,x)(f)|$,
\endroster
then
\roster
\item $\overline{\;\Delta_{A_1}(S\times
 (Z_1\setminus V))\;}\cap T^*(\overline{\;
(\Delta_{A_2}(S\times W))\;})=\emptyset$.
\item $\Delta_{A_1}(S\times (M(A_1)\cap V))
\supset T^*(\Delta_{A_2}(S\times (M_T(A_1)\cap W)))$.
\endroster
\endproclaim
The assertion \therosteritem1 of this proposition
 is a direct consequence of proposition~4.1 when
 we put $Q=Z_2\setminus W$, $\Cal L=\Delta_{A_1}
 (S\times (Z_1\setminus V))$ and we omit its
 proof. Further, \therosteritem2 follows
 immediately from \therosteritem1 and from the
 definition of $M(A_1)$.

\proclaim{Proposition 6.3}Let $Z$ be any compact
 Hausdorff topological space, $A\in\Un(Z;1)$, and
 $1\in A$. Let's put
$$
D=\{\ell|\;\ell\in\partial(A^*)_1,\ell(1)=1\}.
$$
Then $M(A)=\{z|\;z\in Z,\Delta_A(1,z)\in\ex D\}$.
\endproclaim
\demo{Proof of proposition 6.3}Note $D\supset
 \{\Delta_A(1,z)|\;z\in Z\}$ and moreover, $D$
 is an extreme subset of $(A^*)_1$. Then by
 lemma~4.1, $\ex D=D\cap\ext A^*$. Hence
$$\align
M(A)&=\{z|\,z\in Z,\Delta_A(1,z)\in\ext A^*\}\\
 &=\{z|\,z\in Z,\Delta_A(1,z)\in D\cap\ext A^*\}\\
	 &=\{z|\,z\in Z,\Delta_A(1,z)\in\ex D\}.
\endalign
$$
and so, proposition 6.3 is proved.
\enddemo

\remark{Remark 6.2}It follows from both
 proposition 6.3 and \cite{3,~lemma~4.3} that in
 the special case when $Z$ is a Hausdorff compact,
 $A\in\Un(Z;1)$, $1\in A$ our definition of $M(A)$
 is equvalent to the classical definition of
 Choquet boundary (see for instance,
 \cite{3,~page~310}).
\endremark

\remark{Remark 6.3}Note, in the general case
 under consideration in definition 6.1, $M(A)$ is a boundary of $A$ 
 (see section~2 for the definition of boundary of $A$).
If the family $G=\{f\}\subset A\setminus\{0\}$ is 
not placed over  a compact subset of $Z$ (or briefly, 
``$f\in A, f\ne 0$ is not placed over  a compact subset 
of $Z$''), then a well-known immediate  consequence 
of  the Krein-Milman theorem is
$$
\sup_{z\in M(A)}|\Delta_A(1,z)(f)|=\|f\|_A
$$
and we omit the details. Note, it is possible that 
$M(A)\cap \spm(f)=\emptyset$. But, if $f\in A\setminus\{0\}$ is placed over 
some of  compact subsets of $Z$ then we have the following 
assertion --- corollary~6.1. 
\endremark

\proclaim{Corollary 6.1} Let the triple $(Z,\Delta_A,A)$ have 
the property $(*)$, let $Z$ be a Hausdorff topological 
space, and $\Delta_A$ be continuous.
 Suppose, in addition, $\Ch(A)\ne\emptyset$. If $f\in A$ 
 is such that there exists a compact $K\subset Z$ for which
$$
\sup_{z\in Z\setminus K}|\Delta_A(1,z)(f)|<\|f\|_A,
$$
i.e. $f$ is placed over the compact $K$, then
 $$
Y\cap\spm(f)\ne\emptyset, \qquad \forall Y\in\Ch(A).
$$
\endproclaim 
\demo{Proof of corollary 6.1}The one-dimentional
 space $<f>$ spaned over $f$ is such that $(Z,\Delta,<f>)$ has 
 the property $(*)$, where it is setted 
 $\Delta(\lambda,g)=\Delta_A(\lambda,g)$, 
 $\forall \lambda\in S$ and $\forall g\in <f>$. 
 This space has a compact  boundary and by 
 proposition~6.1.\therosteritem1 we obtain $\Ch(<f>)\ne\emptyset$. 
 Moreover,
 a direct computation shows $M(<f>)=\spm(f)$. Further, 
 $<f>\subset A$ and hence by proposition 6.1.\therosteritem2, 
 $M(<f>)\cap Y\ne\emptyset$, so that
 $Y\cap\spm(f)\ne\emptyset$, $\forall Y\in\Ch(A)$.
 Thus corollary~6.1 is proved.
\enddemo

\head 7. The Spaces of Continuous Functions $C(X)$, $C_0(Y)$.\endhead
In this section it is demonstrated an application of the developped above 
theory. Let $X$ be any compact Hausdorff topological space. The  notation 
$C(X)$ stands for the Banach space of all continuous complex valued functions 
on $X$, and this space is normed by the usual  sup-norm. Let $Y$ be any 
locally compact Hausdorff topological space. The notation $C_0(Y)$ stands 
for the Banach space of all continuous complex valued functions on $Y$ 
that vanish at infinity, and this space is normed with the usual sup-norm. 
In this section we obtain all basical results on into-isometries between 
spaces such as $C(X)$, $C_0(Y)$. The theorems stated  below, 
all are simple consequences  of the  theory presented in the previous sections.
 
We define the abstract families $\Cal A$  and $\Cal B$ by
$$\align
\Cal A&=\{(C(X),X)|\;X\text{ is a compact Hausdorff topological space }\},\\
\Cal B&=\{(C_0(Y),Y)|\;Y\text{ is a locally compact Hausdorff topological 
space }\}.
\endalign
$$
If  $(A,Z)\in\Cal A\cup\Cal B$, then the map $\Delta_A:S\times Z\to 
(A^*)_1$ is defined by 
$$
\Delta_A(\lambda,z)(f)=\lambda f(z),
$$
where $\lambda\in S$, $z\in Z$, $f\in A$. Then the triple 
$(Z,\Delta_A,A)$ has  the property $(*)$ and $\Delta_A$ is continuous. 
In particular, $A\in\Un(Z,1)$ accordingly to the definitions of section~2.   

\proclaim{Lemma 7.1}Let $A$ be a  non zero Banach space, 
$(B,Z)\in\Cal A\cup\Cal B$. Suppose that $\Intois(A,B)\ne\emptyset$ 
and let $T\in\Intois(A,B)$ be any. Then for every functional 
$\ell\in\partial (A^*)_1$ there exists a compact $K\subset Z$ such  that 
$$
\ell\not\in T^*(\overline{\;\Delta_B(S\times (Z\setminus K))\;}).
$$
\endproclaim
\demo{Proof of lemma 7.1}Let $\ell\in\partial (A^*)_1$ be any. 
Then, in particular, $\ell\ne  0$ and hence, there is a function $f\in A$ 
such that $\ell(f)\ne 0$. Note $Tf\in B$. So, that 
$$
K=\{z|\;|Tf(z)|\ge |\ell(f)|/2\;\}
$$
is a compact subset of $Z$. Then 
$$
 \sup_{z\in Z\setminus K}|T^*(\Delta_B(1,z))(f)|<|\ell(f)|.
$$
Further, in proposition~4.1 when $A_1=A$, $Z_2=Z$, $A_2=B$, 
$Q=K$, $\Cal L=\{\ell\}$, it is fulfilled the condition of the first bullet. 
Hence by proposition~4.1 it follows  
$$
\{\ell\}\cap T^*(\overline{\;\Delta_B(S\times (Z\setminus K))\;})=\emptyset.
$$
Thus lemma~7.1 is proved.
\enddemo

\proclaim{Lemma 7.2}Let $A$ be a  non zero Banach space, 
$(B,Z)\in\Cal A\cup\Cal B$. Suppose that $\Is(A,B)\ne\emptyset$ 
and let $T\in\Is(A,B)$ be any. Then for every $z_0\in Z$ there exists 
a family of functions $G$, $\emptyset\ne G\subset A\setminus \{0\}$ 
such that  its  image $TG$ is placed over some compact subset of $Z$ 
and in addition 
$$
\bigcap_{f\in G}\spm(Tf)=\{z_0\}.
$$
\endproclaim
\demo{Proof of lemma 7.2}There are two cases: $Z=\{z_0\}$, and 
$Z\setminus\{z_0\}\ne\emptyset$. In the first case lemma~7.2 is 
a trivial one and we omit the details of the proof. So, let's consider 
the second case. Let $z_0\in  Z$ be any. Then there  exists an 
open neighbourhood $U$ of $z_0$ such that its closure $\overline U$ is  
compact. Further, it follows by the Urysohn's lemma that there is a  
function $h_0\in B$ such that  $z_0\in\spm(h_0)$, 
$\sup_{z\in Z\setminus\overline U}|\Delta_B(1,z)(h_0)|=0$. 
Moreover, by the Urysohn's lemma, there exists a  family of functions 
$H\subset B$ such that $\bigcap_{h\in H}\spm(h)=\{z_0\}$. Thus 
the family $\{h_0\}\cup H$ is placed over the compact $\overline U$ 
and is such that
$$
\bigcap_{h\in \{h_0\}\cup H}\spm(h)=\{z_0\}.
$$
Hence $G=\{f|\;f\in A,Tf\in \{h_0\}\cup H\}$ has the desired properties.

Thus lemma~7.2 is proved.        
\enddemo

Recall that a map $F:\Omega_1\to\Omega_2$ is said to be perfect \cite{6,~\S~3.7}  if and  only if the following three conditions are fulfilled:
\roster
\item"{$\bullet$}" $\Omega_1$ is a Hausdorff topological space and
 $\Omega_2$ is a topological space;
\item"{$\bullet$}" $F$ is a contionuous closed map;
\item"{$\bullet$}" the whole preimage
 $F^{-1}(\{y\})=\{x|\;x\in \Omega_1,\;F(x)=y\}$ is a compact subset of $ \Omega_1$, for every point $y$ in  the image $F(X)=\{F(x)|\;x\in \Omega_1\}$ of  the map $F$.  
\endroster 

\proclaim{Lemma 7.3. A topological lemma}Let for $i=1,2$ the set $\Omega_i$ be a Hausdorff topological space and let $F:\Omega_1\to\Omega_2$ be a surjective continuous map. If for every $y\in\Omega_2$ there is a compact $K\subset \Omega_1$ such that $y\not\in F(\Omega_1\setminus K)$, then the map $F$ is perfect.
\endproclaim
\demo{Proof of lemma 7.3}We prove this lemma directly by using the  definition of perfect  map. 

First,  we prove that $F$ is closed. Let the closed set $W\subset \Omega_1$ be any,   
and let $z$ be any point of the  closure $\overline{\;F(W)\;}$ of $F(W)$. By the assumption there  is a compact $K\subset \Omega_1$ such that $z\not\in F(\Omega_1\setminus K)$. Hence,
$$\align
z\in \overline{\;F(W)\;}&\subset \overline{\;F(W\cap K)\;}\cup 
\overline{\;F(W\setminus K)\;}\\
                                          &\subset \overline{\;F(W\cap K)\;}\cup
\overline{\;F(\Omega_1\setminus K)\;},
\endalign
$$ 
so, that $z\in \overline{\;F(W\cap K)\;}$. Further, it follows by continuity of  $F$ and by  compactness of $W\cap K$, that $\overline{\;F(W\cap K)\;}=F(W\cap K)$. Therefore, $z\in F(W\cap K)$ and in particular, $z\in F(W)$. Accordingly to the choise of $z\in \overline{\;F(W)\;}$, it follows that  $\overline{\;F(W)\;}\subset F(W)$. So, we obtain that the map  $F$ is  closed.

Second, we prove that the preimage $F^{-1}(y)$ is compact, $\forall y\in \Omega_2$. Let $y\in \Omega_2$ be  any. Then from surjectivity of $F$ it follows  $F^{-1}(y)\ne\emptyset$. Further, by the assumption, there is a  compact $K\subset \Omega_1$ such that $y\not\in F(\Omega_1\setminus K)$. Then $F^{-1}(y)\cap (\Omega_1\setminus K)=\emptyset$. So, $F^{-1}(y)\subset K$. Further, by continuity of $F$ it follows that $F^{-1}(y)$ is closed. Hence $F^{-1}(y)$ is compact. 

Then by the definition, we  obtain that $F$ is perfect. Thus lemma~7.3 is proved.      
\enddemo 

\proclaim{Proposition 7.1}Let $A$ be a non zero Banach
 space, the triple $(Z,\Delta_B,B)$
 have the property $(*)$ and let $\Delta_B$ be a
 continuous map. If the isometry $T\in\Intois(A,B)$ is such that  
for  every $\ell\in\partial (A^*)_1$ there is a compact $K\subset Z$ 
so that the following condition is satisfied
$$
\ell\not\in T^*(\overline{\;\Delta_B(S\times (Z\setminus K))\;})\tag{7.1}
$$  
then $\ext A^*\subset T^*\Delta_B(S\times Z)$.
\endproclaim
\demo{Proof of proposition 7.1}Let $\ell\in\ext A^*$ be any. 
Then by theorem~3.1 it follows 
$$
\ell\in T^*(\overline{\;\Delta_B(S\times Z)\;})
$$ 
Further, by the well-known inclusion $\ext  A^* \subset \partial 
(A^*)_1$ we  obtain $\ell\in \partial (A^*)_1$. Hence there is a  
compact $K\subset Z$ such  that the condition \thetag{7.1} is  fulfilled. Therefore,  
$$
\ell\in T^*(\overline{\;\Delta_B(S\times K)\;}).
$$ 
Note $\overline{\;\Delta_B(S\times K)\;}=\Delta_B(S\times K)$
 because of compactness of $K$ and of continuity of  $\Delta_B$. So, that 
$$
\ell\in T^*\Delta_B(S\times K).
$$ 
Hence, accordingly to the choise  of $\ell$, we obtain 
$\ext A^*\subset T^*\Delta_B(S\times Z)$. 

Thus proposition~7.1 is proved.  
\enddemo

\proclaim{Proposition 7.2}Let $A$ be a non zero Banach
 space, the triple $(Z,\Delta_B,B)$
 have the property $(*)$ and let $\Delta_B$ be a
 continuous map. If the isometry 
$T\in\Intois(A,B)$ is such that for every $z_0\in Z$ there is a family of 
functions $G$, with the following properties $\emptyset\ne 
G\subset  A\setminus\{0\}$, its image $TG$ is placed over 
some compact subset of $Z$, and 
$$
\bigcap_{f\in G}\spm(Tf)=\{z_0\},
$$
then the map $T^*\Delta_B:S\times Z\to \ext A^*$ is well defined, 
injective, and continuous.
\endproclaim
We sketch the main arguments and omit the details. By 
corollary~4.1 it follows that this map is well defined. 
Further, by the assumptions on the isometry $T$ it follows that 
$|TA|$ distinguishes all points of $Z$ (see definition~5.2). 
Hence by proposition~5.2 it follows that the map $T^*\Delta_B$ 
is injective. And of course, it is continuous because 
$T^*$ and $\Delta_B$, both are continuous. These arguments are 
enough for the proof of  proposition~7.2.

\definition{Definition 7.1}Let $A$ be a non zero Banach
 space, the triple $(Z,\Delta_B,B)$
 have the property $(*)$ and let $\Delta_B$ be a
 continuous map. Then we define the following two properties $(\alpha)$  and $(\beta)$  of an isometry $T\in \Intois(A,B)$:
\roster
\item"{$(\alpha)$}" for  every $\ell\in\partial (A^*)_1$ there is a compact $K\subset Z$ 
so that the following condition is satisfied
$$
\ell\not\in T^*(\overline{\;\Delta_B(S\times (Z\setminus K))\;})
$$       
\item"{$(\beta)$}"  for every $z_0\in Z$ there is a family of 
functions $G$, with the following properties $\emptyset\ne 
G\subset  A\setminus\{0\}$, its image $TG$ is placed over 
some compact subset of $Z$, and 
$$
\bigcap_{f\in G}\spm(Tf)=\{z_0\},
$$
\endroster
\enddefinition

\proclaim{Theorem 7.1}Let $A$ be a non zero Banach space, 
the triple $(Z,\Delta_B,B)$ have the property $(*)$ and let $\Delta_B$ be a
 continuous map. Then we have the following assertions:
\roster
\item if $T\in\Intois(A,B)$ has the  property $(\alpha)$ then $\Ch_T(A)\ne\emptyset$ and the map 
$$
T^*\Delta_B:S\times M_T(A)\to \ext  A^*
$$
is well  defined, surjective and perfect;
\item if $T\in\Intois(A,B)$ has the two properties $(\alpha)$, $(\beta)$ then 
$M_T(A)=Z$ and the map 
$$
T^*\Delta_B:S\times Z\to \ext  A^*
$$
is a homeomorphic map.
\endroster 
\endproclaim
\demo{Proof of  theorem  7.1}The assertion of~\therosteritem1 follows immediately from lem\-ma~7.3 and from proposition~7.1. The assertion of~\therosteritem2 is a  direct consequence of~\therosteritem1 and of  proposition~7.2.

Thus theorem~7.1 is poved.  
\enddemo

\proclaim{Corollary 7.1}Let for $i=1,2$, the couple   $(A_i,Z_i)\in\Cal A\cup\Cal B$. Then we have the following assertions:
\roster
\item if $T\in\Intois(A_1,A_2)$ then $\Ch_T(A)\ne\emptyset$ and the map 
$$
T^*\Delta_{A_2}:S\times M_T(A_1)\to \ext  A_1^*
$$
is well  defined, surjective and perfect;
\item if $T\in\Is(A_1,A_2)$ then 
$M_T(A_1)=Z_2$ and the map 
$$
T^*\Delta_{A_2}:S\times Z_2\to \ext  A_1^*
$$
is a homeomorphic map.
\endroster 
\endproclaim
\demo{Proof of corollary 7.1}Note, that for $i=1,2$, the map $\Delta_{A_i}:S\times Z_i\to (A_i^*)_1$ is  defined by $\Delta_{A_i}(\lambda,x)(f)=\lambda f(x)$, $\forall (\lambda,x)\in S\times Z_i$ and $\forall f\in A_i$. 

First we prove \therosteritem1. We apply lemma~7.1 in the case $Z=Z_2$, $B=A_2$, $A=A_1$, $\Delta_B=\Delta_{A_2}$. Thus we obtain that the isometry $T\in \Intois(A_1, A_2)$ has the property $(\alpha)$. Then by  theorem~7.1.\therosteritem1, with $Z=Z_2$, $B=A_2$, $A=A_1$, $\Delta_B=\Delta_{A_2}$, it follows that the map $T^*\Delta_{A_2}:S\times Z_2\to \ext A_1^*$ is well defined, surjective and  perfect. So, \therosteritem1 is proved. 

Second, we prove \therosteritem2.  We apply both lemma~7.1 and lemma~7.2  in the case $Z=Z_2$, $B=A_2$, $A=A_1$, $\Delta_B=\Delta_{A_2}$. Thus we obtain that the isometry $T\in \Is(A_1, A_2)$ has the properties $(\alpha)$, $(\beta)$. Then by  theorem~7.1.\therosteritem2, with $Z=Z_2$, $B=A_2$, $A=A_1$, $\Delta_B=\Delta_{A_2}$, it follows that the map $T^*\Delta_{A_2}:S\times Z_2\to \ext A_1^*$ is homeomorphic. Thus \therosteritem2 is proved and  the corollary~7.1 is proved, too.  
\enddemo

\proclaim{Theorem 7.2}Let for $i=1,2$  $(A_i,Z_i)\in\Cal A\cup\Cal B$. Suppose $\Intois(A_1,A_2)\ne\emptyset$. Then  for every $T\in\Intois(A_1,A_2)$ there exists a  couple $(\phi,\tau)$ of maps such that:
\roster
\item $\phi:M_T(A_1)\to S$ is continuous;
\item $\tau:M_T(A_1)\to Z_1$ is surjective perfect map;
\item $Tf(x)=\phi(x)f(\tau(x)),\quad\forall f\in A_1, \forall x\in M_T(A_1)$;
\endroster
 This couple is unique in the following sence: if the couple $(\phi',\tau')$ is such that 
$$\align
&\phi':U\to S,\\
&\tau':U\to Z_1,\\
&Tf(x)=\phi'(x)f(\tau'(x)),\qquad\forall x\in U, \forall f\in A_1
\endalign
$$ 
where $\emptyset\ne U\subset Z_2$, then $U\subset M_T(A_1)$ and $\phi'(x)=\phi(x)$, $\tau'(x)=\tau(x)$, for all $x\in U$.  
\endproclaim

\remark{Remark}We have to recall that in consequence of theorem~7.2.\therosteritem2 we obtain some of topological properties of $M_T(A_1)$. We use the topological results of \cite{6,~\S 3.7} to obtain the following:
\roster
\item $M_T(A_1)$ is compact $\iff$ $Z_1$ is compact;
\item $M_T(A_1)$ is locally compact (but not compact) $\iff$ $Z_1$ is locally compact (but not compact, resp.);
\item if $M_T(A_1)$ is not compact but is closed subset of $Z_2$, then there exists a net $\{x_\alpha\}$ of points of  $M_T(A_1)$ that tends  to the infinity $\infty_2$ of $Z_2$. Further, by \therosteritem2, $Z_1$ is not compact. Let $\infty_1$ stand for the infinity of $Z_1$. Then we have
$$
\tau(x)\longrightarrow\infty_1 \Longleftrightarrow x\longrightarrow\infty_2.
$$
\endroster
\endremark

\remark{Remark}In the proof of theorem~7.2 we make use of the following projections:
$$\matrix \format \r&\l&\quad\r&\l&\quad\l\\
q&:S\times Z_1\to S,&q(\lambda,z)&=\lambda,&\forall (\lambda,z)\in S\times Z_1;\\
p_1&:S\times Z_1\to Z_1,&p_1(\lambda,z)&=z,&\forall (\lambda,z)\in S\times Z_1;\\
p_2&:S\times M_T(A_1)\to M_T(A_1),&p_2(\mu,y)&=y,&\forall (\mu,y)\in S\times M_T(A_1),
\endmatrix
$$
where for $i=1,2$  $(A_i,Z_i)\in\Cal A\cup\Cal B$, $T\in\Intois(A_1,A_2)$.
In particular, $p_1$, $p_2$, both are surjective perfect  maps and $q$ is surjective continuous quotient map. Note that $q$  is perfect if and only if $Z_1$ is compact Hausdorff topological space.
\endremark

\demo{Proof of theorem 7.2}Let $\ext A_1^*$ has the induced weak${}^*$-topology of $A_1^*$. By corollary~7.1.\therosteritem1 we obtain $\Ch_T(A_1)\ne\emptyset$ and 
$$
T^*\Delta_{A_2}:S\times M_T(A_1)\to \ext A_1^*
$$
is well defined, surjective and perfect. Further, by corollary~7.1.\therosteritem2 
$$
\Delta_{A_1}:S\times Z_1\to \ext A_1^*
$$ 
is a homeomorphic map. Therefore the map
$$
\Delta_{A_1}^{-1}T^*\Delta_{A_2}:S\times M_T(A_1)\to \ext A_1^*
$$
is surjective and perfect.

We set
$$\align
\phi_0&=q\Delta_{A_1}^{-1}T^*\Delta_{A_2},\\
\tau_0&=p_1\Delta_{A_1}^{-1}T^*\Delta_{A_2}.
\endalign
$$
Accordingly, $\tau_0$ is surjective and  perfect. Further,
$$
\align
(\phi_0(\lambda,x),\tau_0(\lambda,x))&=\Delta_{A_1}^{-1}T^*\Delta_{A_2} (\lambda,x)=\Delta_{A_1}^{-1}T^*(\lambda\Delta_{A_2}(1,x))\\
&=\Delta_{A_1}^{-1}(\lambda T^*\Delta_{A_2}(1,x))=
\Delta_{A_1}^{-1}(\lambda(\Delta_{A_1}\Delta_{A_1}^{-1}T^*\Delta_{A_2}
(1,x)))\\
&=\Delta_{A_1}^{-1}(\lambda(\Delta_{A_1}(\phi_0(1,x),\tau_0(1,x))))\\
&=\Delta_{A_1}^{-1}(\Delta_{A_1}(\lambda\phi_0(1,x),\tau_0(1,x)))\\
&=(\lambda\phi_0(1,x),\tau_0(1,x)),\qquad \forall(\lambda,x)\in S\times M_T(A_1).
\endalign
$$
We define the maps $\phi:M_T(A_1)\to S$, $\tau:M_T(A_1)\to Z_1$ by
$$
\align
\phi(x)&=\phi_0(1,x),\\
\tau(x)&=\tau_0(1,x),
\endalign
$$
where $x\in M_T(A_1)$. Then, $\tau$ is surjective and perfect, $\phi$ is continuous. We have to  remark that $\phi$ is perfect iff $Z_1$ is compact (in the case when $Z_1$ is compact then the diagonal map $(\phi,\tau):M_T(A_1)\to S\times Z_1$ is perfect and is not necessary surjective). 

Further, according to the definition of the maps  $\phi$, $\tau$
$$
\Delta_{A_1}^{-1}T^*\Delta_{A_2}(1,x)=(\phi(x),\tau(x)), \qquad\forall x\in M_T(A_1).
$$ 
Then $T^*\Delta_{A_2}(1,x)=\Delta_{A_1}(\phi(x),\tau(x))$ so, that
$$
Tf(X)=\phi(x)f(\tau(x)),\qquad \forall x\in M_T(A_1), \forall f\in A_1.
$$ 
Thus it is proved the existence of maps $\phi$, $\tau$ for which the assertions  \therosteritem1, \therosteritem2, \therosteritem3 of theorem~7.1 hold on. 
We claim that such a couple $(\phi,\tau)$ of maps is unique for the isometry $T$. Indeed, for any couple $(\phi',\tau')$ with the properties:
$$\align
&\phi':U\to S,\\
&\tau':U\to Z_1,\\
&Tf(x)=\phi'(\tau'(x)),\qquad\forall x\in U, \forall f\in A_1,
\endalign 
$$
we have 
$$
T^*\Delta_{A_2}(1,x)=\Delta_{A_1}(\phi'(x),\tau'(x)),\qquad x\in U.
$$
Note, by corollary~7.1 $\Delta_{A_1}(\phi'(x),\tau'(x))\in\ext  A_1^*$, $x\in U$. So, $T^*\Delta_{A_2}(1,x)\in\ext  A_1^*$ and  hence $x\in M_T(A_1)$. Thus $U\subset M_T(A_1^*)$. Further
$$
(\phi'(x),\tau'(x))=\Delta_{A_1}^{-1}T^*\Delta_{A_2}(1,x)=(\phi(x),\tau(x)),  \qquad \forall x\in U.
$$
Therefore the uniqueness is proved. 

So, theorem~7.2 is proved.
\enddemo

\proclaim{Lemma 7.4}Let for $i=1,2$  $(A_i,Z_i)\in\Cal A\cup\Cal B$. Suppose there is  a non closed subset $U$ of $Z_2$ and there is a surjective perfect map $\pi: U\to Z_1$. Then for  every $x_0\in \overline U\setminus U$ and for every $f\in A_1$ 
$$
\lim_{x\to x_0, x\in U}f(\pi(x))=0.
$$
where $\overline U$ stands for  the  closure of $U$ in $Z_2$.
\endproclaim
\demo{Proof of lemma 7.4}Let $x_0\in \overline U\setminus U$, $f\in A_1$,  $\varepsilon>0$ all be choosen in an arbitrary manner. Then 
$$
K=\{y|\;y\in Z_1,\;|\Delta_{A_1}(1,x)(f)|\ge \varepsilon\}
$$
is compact. Hence its preimage $\pi^{-1}(K)$ is a compact subset of $U$, because of the perfectness of $\pi$. So, that $\pi^{-1}(K)$ is closed in $\overline U$. Moreover, it follows from $x_0\not\in U$ that  $x_0\not\in \pi^{-1}(K)$. Therefore $\overline U\setminus \pi^{-1}(K)$ is an open neighbourhood of  $x_0$ in $\overline U$ with  respect to the subspace topology induced  by $Z_2$. Then 
$$
\lim_{x\to x_0, x\in  U}|f(\pi(x))|\le\sup_{x\in U\setminus \pi^{-1}(K)}|f(\pi(x))|=\sup_{z\in Z_1\setminus K}|f(z)|\le\varepsilon.
$$
Thus $\displaystyle\lim_{x\to x_0, x\in  U}|f(\pi(x))|=0$ and hence lemma~7.4 is proved.
\enddemo

\proclaim{Corollary 7.2}Let for $i=1,2$ the couple  $(A_i,Z_i)\in\Cal A\cup\Cal B$. Suppose $\Intois(A_1,A_2)\ne\emptyset$. If $T\in\Intois(A_1,A_2)$ is such that $M_T(A_1)$ is not closed then for every $f\in A_1$
$$
Tf(x)=0,\qquad\forall x\in \overline{M_T(A_1)}\setminus M_T(A_1).
$$
\endproclaim
Instead of an explicite proof we note that this corollary is an immediate consequence of both theorem~7.2 and lemma~7.4. 

\proclaim{Corollary 7.3}Let for $i=1,2$  $(A_i,Z_i)\in\Cal A\cup\Cal B$. If there exist both an open subset $U$ of $Z_2$ and a surjective perfect map $\pi:U\to Z_1$, then $\Intois(A_1,A_2)\ne\emptyset$.
\endproclaim
Instead of a detailed proof we remark  the following arguments. By lemma~7.4, for every $f\in A_1$ the function 
$$
(\pi^of)(x)=\cases
                   f(\pi(x)),&x\in U\\
                   0, &x\not\in U
\endcases
$$ 
is continuous on $Z_2$. Hence the linear operator $\pi^0:f\mapsto \pi^0f$ belongs  to $\Intois(A_1,A_2)$. These arguments prove corollary~7.3. 

\remark{Remark 7.1}The problem when $\Intois(A_1,A_2)\ne\emptyset$ under consideration in corollary~7.3 is not new. There are results in \cite{10} when $U$ is closed.
\endremark

\proclaim{Corollary 7.4}Let for $i=1,2,3$ the couple  $(A_i,Z_i)\in\Cal A\cup\Cal B$. Suppose $\Intois(A_1,A_2)\ne\emptyset$, $\Intois(A_2,A_3)\ne\emptyset$. Let  $T_1\in\Intois(A_1,A_2)$, $T_2\in\Intois(A_2,A_3)$ both be  any and let's put $T_3=T_2T_1$ (so, that  $T_3\in\Intois(A_1,A_3)$). Let for $i=1,2,3$, the couple of maps $(\phi_i,\tau_i)$ be the corresponding couple to $T_i$ which existence is established in theorem~7.2. Then:
\roster
\item $\emptyset\ne \tau_2^{-1}(M_{T_1}(A_1))\subset M_{T_2}(A_2)\cap M_{T_3}(A_1)$; 
\item $\tau_2^{-1}(M_{T_1}(A_1))\in \Ch_{T_3}(A_1)$, and $M_{T_2}(A_2)\cap M_{T_3}(A_1)\in \Ch_{T_3}(A_1)$;
\item $\phi_3(x)=\phi_2(x)\phi_1(\tau_2(x))$, for $x\in \tau_2^{-1}(M_{T_1}(A_1))$;
\item $\tau_3(x)=\tau_1(\tau_2(x))$, for $x\in \tau_2^{-1}(M_{T_1}(A_1))$;
\item $\tau_3:\tau_2^{-1}(M_{T_1}(A_1))\to Z_1$ is surjective and perfect;
\endroster  
\endproclaim
\demo{Proof of corollary 7.4}By theorem~7.2 
\roster
\item"{$\bullet$}" $\tau_1:M_{T_1}(A_1)\to  Z_1$, $\tau_1$ is  surjective and perfect; 
\item"{$\bullet$}" $\tau_2:M_{T_2}(A_2)\to  Z_2$, $\tau_2$ is  surjective and perfect;
\item"{$\bullet$}" $\tau_3:M_{T_3}(A_1)\to  Z_1$, $\tau_3$ is  surjective,
\endroster 
where $M_{T_1}(A_1)\subset Z_2$, $M_{T_2}(A_2)\subset Z_3$, and $M_{T_3}(A_1)\subset Z_3$. In particular, $\emptyset\ne\tau_2^{-1}(M_{T_1}(A_1))\subset M_{T_2}(A_2)$, and by the definition $T_3=T_2T_1$ it follows
$$\align
(T_3f)(x)&=(T_2(T_1f))(x)=\phi_2(x)(T_1f)(\tau_2(x))\\
               &=\phi_2(x)\phi_1(\tau_2(x))f(\tau_1(\tau_2(x))),
\endalign
$$
where $x\in \tau_2^{-1}(M_{T_1}(A_1))$, and $f\in A_1$ both are arbitrary. Then by theorem~7.2 (the part about uniquness) we  obtain $\tau_2^{-1}(M_{T_1}(A_1))\subset  M_{T_3}(A_1)$ and hence corollary~7.4.\therosteritem1 is proved. Note, in addition, the assertions of   corollary~7.4.\therosteritem3 and of corollary~7.4.\therosteritem4, are proved,  too.  Further, from \therosteritem4 we obtain the part  of \therosteritem5 on surjectivity. The assertion on perfectness we obtain by using the general properties of perfect maps \cite{6,~\S~3.7}: the restriction of the  perfect map $\tau_2:M_{T_2}(A_2)\to Z_2$ on the  whole preimage $\tau_2:\tau_2^{-1}(M_{T_1}(A_1))\to M_{T_1}(A_1)$ is perfect, too; the composition $\tau_3=\tau_1\circ\tau_2$ of perfect maps is perfect, too. So, \therosteritem5 is proved.  Note, in particular, 
$\tau_3:\tau_2^{-1}(M_{T_1}(A_1))\to Z_1$ is surjective. Recall, by corollary~7.1.\therosteritem2, we have 
$\Delta_{A_1}:S\times Z_1\to \ext A_1^*$ is bijective. Therefore,     
$\tau_2^{-1}(M_{T_1}(A_1))\in \Ch_{T_3}(A_1)$ and so, the first part of \therosteritem2 is proved. The second one follows from both, from this one and from \therosteritem1. Thus corollary~7.4 is proved.
\enddemo

\proclaim{Corollary 7.5}Let for $i=1,2$ the couple  $(A_i,Z_i)\in\Cal A\cup\Cal B$. Suppose $\Is(A_1,A_2)\ne\emptyset$ and $T\in \Is(A_1,A_2)$ be any. Then 
\roster
\item $\tau_T:Z_2\to Z_1$ is a homeomorphic map;
\item $\tau_{T^{-1}}=\tau_T^{-1}$, $\phi_{T^{-1}}=\dfrac1{\phi_T\circ\tau_{T^{-1}}}$,
\endroster
where the two couples of maps $(\phi_T,\tau_T)$ and  $(\phi_{T^{-1}},\tau_{T^{-1}})$ are the corresponding couples to the isometries $T$ and $T^{-1}$, resp. (see theorem~7.2).    
\endproclaim
Instead of a detailed proof we note that corollary~7.5 is a direct consequence of   corollary~7.4.

\enddocument